\documentclass[11pt,a4paper,twoside]{article}
\usepackage{a4wide}
\usepackage{amscd,amssymb,amsmath,amsthm,latexsym}
\usepackage[all]{xy}
\usepackage[mathscr]{euscript}
\usepackage{epsfig}
\usepackage[dvips]{color}

\newcommand{\oGamma}{\overline \Gamma}

\DeclareMathOperator{\mTr}{Tr}

\newcommand{\str}{\stackrel}
\newcommand{\lra}{\longrightarrow}

\newcommand{\Ho}{\mbox{\rm H}}

\newcommand{\mZ}{\mathcal C}
\newcommand{\mQ}{\mathcal Q}

\newcommand{\Ext}{\mbox{\rm Ext}_{H^e}^1}

\newcommand{\hin}{\makebox(0,0)}

\newcommand{\Rr}{{\mathbb R}}

\newcommand{\Zz}{{\mathbb Z}}

\newcommand{\mand}{\mbox{ and }}

\newcommand{\Aa}{{\Bbb A}}
\newcommand{\Cc}{{\Bbb C}}

\newcommand{\Pp}{{\Bbb P}}

\newtheorem{Prop}{Proposition}[section]
\newtheorem{Thm}[Prop]{Theorem}
\newtheorem{Lemma}[Prop]{Lemma}
\newtheorem{Cor}[Prop]{Corollary}

\begin{document}
\title{Cluster-Cyclic Quivers with three Vertices and the Markov Equation, \\~ \\with an appendix by Otto Kerner}
\author{Andre Beineke \and Thomas Br\"ustle \and Lutz Hille}
\date{\today}
\maketitle
\begin{abstract}
Acyclic cluster algebras have an interpretation in terms
 of tilting objects in a Calabi-Yau category defined by some hereditary algebra. For a given quiver $Q$ it is thus desirable to decide if the cluster algebra defined by $Q$ is acyclic. 
We call  $Q$ 
cluster-acyclic in this case, otherwise cluster-cyclic.
  In this note we classify the
cluster-cyclic quivers with three vertices using a Diophantine equation studied by Markov.
\end{abstract}

\section{Introduction}\label{sectintro}
Cluster algebras have been introduced  and studied by Fomin and
Zelevinsky in \cite{FZ1,FZ2}.  
In \cite{BMRRT} it is shown that each
acyclic cluster algebra admits an interpretation in terms
 of tilting objects in a so-called cluster category. This category  is a Calabi-Yau category defined as a quotient of the derived category of modules over some hereditary algebra. 
It is thus desirable to study which cluster algebras are acyclic and which are not. While it has been shown in \cite{FZ2} that all finite cluster algebras are acyclic, the general case is not known.

The first non-trivial situation is that of cluster algebras of rank three. 
This case has recently been studied in \cite{K} (which analyzes properties
of the acyclic case by represen\-tation-theoretic methods) and in
\cite{ABBS} (which describes the mutation classes of non acyclic cluster
algebras). The cluster algebras considered here are defined by a
quiver $Q$; the quiver $Q$ corresponds to a skew-symmetric matrix
$B$. In general, a cluster algebra is given by a skew-symmetrizable
matrix. 
In this paper we give precise criteria for deciding which quivers with
three vertices yield an acyclic cluster algebra (Theorem \ref{Tmain1})
and which do not (Theorem \ref{Tmain2}). 
  
We consider a quiver $Q$ with three vertices $1,2$ and $3$ and 
cyclic orientation. 
We assume $Q$ has $x$ arrows from $1$ to $2$,
$y$ arrows from $2$ to $3$ and $z$ arrows from $3$ to $1$:

$$  
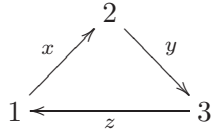
\xymatrix{ & 2 \ar^{y}[dr] & \\ 1  \ar^{x}[ur]  & & 3 \ar^{z}[ll]  }  $$

An essential tool in the definition of a cluster algebra is the
mutation of a quiver.  
We describe the effect of these cluster mutations on the cyclic  quiver $Q$.
A cluster mutation in the vertex $2$ defines a new quiver $\mu_2 Q$
which is obtained from $Q$ by reversing all arrows from $1$ to $2$ and
from $2$ to $3$. 
Moreover, the quiver $\mu_2Q$ has  $z-xy$ arrows from $3$  
to $1$, provided this number is non-negative. Then the quiver $\mu_2Q$
is no longer 
cyclic, it is {\sl acyclic}. If $z-xy$ is negative, then
$\mu_2Q$ has $xy-z$ arrows from $1$ to $3$ and it is cyclic:

$$
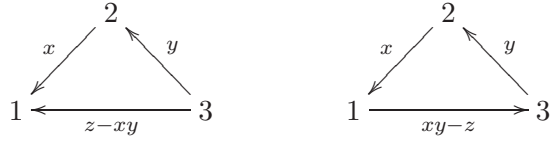
\xymatrix{ & 2 \ar_{x}[dl] & \\ 1 & & 3  \ar_{y}[ul]
  \ar^{z-xy}[ll] }  \qquad \qquad 
 \xymatrix{ & 2 \ar_{x}[dl] & \\ 1  \ar_{xy-z}[rr] & & 3  \ar_{y}[ul]  } $$

The cluster mutations $\mu_1Q$ and $\mu_3Q$ in the vertices $1$ and
$3$ are defined analogously.  
We say that $Q$ is {\sl cluster-acyclic} if
there exists  a finite sequence of cluster
mutations from $Q$ to an acyclic quiver $Q'$. Otherwise, if all
sequences of cluster mutations 
yield cyclic quivers, we say that $Q$ is {\sl cluster-cyclic}. 
Of course, these cluster mutations are compatible
with  permutations of the triple $(x,y,z)$ defining the cyclic quiver
$Q$. So we say that a triple of non-negative integers $(x,y,z)$ is
{\sl cluster-acyclic} (respectively, {\sl cluster-cyclic}) if the
corresponding quiver $Q$ has this property. 

In this note, we  characterize cluster-cyclic quivers in terms of a
Diophantine equation first studied 
by  Markov in \cite{Markov}. 
For
each triple $(x,y,z)$, we define its {\sl Markov constant} as
$$
C(x,y,z) := x^2 + y^2 + z^2 -xyz.
$$
It turns out that, with some exceptions, the value of the Markov
constant characterizes cluster-acyclic
quivers. We first note that $C(x,y,z)$ is invariant under cluster mutations
(Lemma \ref{LMarkovinv}). We  consider the action of the group
$\oGamma$  on $\Rr^3$, where
$\oGamma$ is  generated by the three cluster mutations $\mu_1, \mu_2$
and $\mu_3$ and all permutations
in $S_3$. For each $C \in \Rr$, this group action reduces to an action on the affine algebraic
variety 
$$
V(C) := \{ (x,y,z) \in \Aa^3 \mid C(x,y,z) = C \}.
$$
We are mainly interested in the integral points of $V(C)$, but will
also consider it over $\Rr$ and over
$\Cc$.
Moreover, for $x,y \geq 2$, we  define two functions
$$
m^-(x,y) := \frac{1}{2} (xy - \sqrt{(x^2 - 4)(y^2 - 4)}) \mand m^+(x,y) := 
\frac{1}{2} (xy + \sqrt{(x^2 - 4)(y^2 - 4)}),
$$
where $m^+(x,y) \geq m^-(x,y) \geq 2$ for all $x,y \geq 2$.
If we consider a point in $V(4)$  and express $z$ as a
function of  $x$ and $y$, then we
obtain the two solutions $z = m^+(x,y)$ and $z = m^-(x,y)$. Moreover,
it is easy to check (using Lemma \ref{LMarkovinv}) 
that
$$
m^+(x,m^-(x,y)) = m^+(x,xy - m^+(x,y)) = y = m^-(x,m^+(x,y)) = m^-(x,
xy - m^-(x,y)).
$$
Using these two functions and the Markov constant, we can characterize
cluster-acyclic triples as follows.
If the Markov constant is larger than $4$, we have a cluster-acyclic
triple; if the Markov constant is
less then $0$, we have an  cluster-cyclic triple. If the Markov constant
lies in the interval $[0,4]$ then
we have a finite number of cluster-acyclic triples for the Markov constants
$0,1,2,4$ corresponding to four finite
orbits (see Theorem \ref{Tmain1} (3) below) and
an infinite number of cluster-cyclic orbits (one contains one element,
all others are infinite
orbits) for  the Markov constants $0$ and $4$ (see Theorem \ref{Tmain2}).
Also note that $4$ is the only value of the Markov constant with an
infinite number of integral orbits; for all other values
the number of orbits is finite (see Corollary \ref{Cfinnumborbits}).

\begin{Thm} \label{Tmain1} 
Let $Q$ be a cyclic quiver with numbers of arrows given by $x,y,z$ in
$\Zz_{\geq 0}^3$. Then the following
conditions are equivalent. \\
(1) The triple $(x,y,z)$ is cluster-acyclic. \\
(2) The Markov constant satisfies $C(x,y,z) > 4$  or $\min\{x,y,z\} < 2$. \\
(3) The Markov constant satisfies $C(x,y,z) > 4$ or the triple $(x,y,z)$
is in the
 following list (where we assume $x \geq y \geq z$): \\
\hspace*{0.3cm} a) $C(x,y,z) = 0:$ $(x,y,z) = (0,0,0)$, \\
\hspace*{0.3cm} b) $C(x,y,z) = 1:$ $(x,y,z) = (1,0,0)$, \\
\hspace*{0.3cm} c) $C(x,y,z) = 2:$ $(x,y,z) = (1,1,0)$ or $(1,1,1)$, \\
\hspace*{0.3cm} d) $C(x,y,z) = 4:$ $(x,y,z) = (2,0,0)$ or $(2,1,1)$. \\
(4) $\min\{x,y,z\} < 2$, or $\min\{x,y \}\geq 2$ and precisely one of
the following inequalities is satisfied: $z > m^+(x,y)$ or $z <
m^-(x,y)$. 
\end{Thm}

\begin{Thm} \label{Tmain2}
Let $Q$ be a cyclic quiver with numbers of arrows given by $x,y,z$ in
$\Zz_{\geq 0}^3$.  
Then the following
conditions are equivalent. \\
(1) The triple $(x,y,z)$ is cluster-cyclic. \\
(2) $x,y,z \geq 2$ and $C(x,y,z) \leq 4$. \\
(3) $C(x,y,z) < 0$ or the triple is in one of the followong orbits: \\
\hspace*{0.3cm} a) $C(x,y,z) = 0:$ $(x,y,z)$ is in $\oGamma (3,3,3)$ or \\
\hspace*{0.3cm} b) $C(x,y,z) = 4:$ $(x,y,z)$ is in $\oGamma (2,x,x)$ for some $x \geq 2$. \\
(4) $x,y,z \geq 2$ and $m^+(x,y) \geq z \geq m^-(x,y)$. \\
(5) $x,y,z \geq 2$ and in the orbit $\oGamma (x,y,z)$ there exists a unique element
$ f(x,y,z) := (x',y',z')$ satisfying $x' \geq y' \geq z'$ and $y'z' - x'
\geq x'$. 
\end{Thm}

\begin{Cor}
(1) The list in Theorem \ref{Tmain1}, (3) together with the orbit of
$(x,0,0)$ (for $x > 2$) is the list of all finite cluster-acyclic
orbits containing a triple with $x\geq y \geq z \geq 0$.  \\
(2) The only finite cluster-cyclic orbit is $(2,2,2)$. 
\end{Cor}

The existence of two functions $m^+$ and $m^-$ satisfying condition
(4) in Theorem \ref{Tmain2} has already been shown in
\cite{ABBS} without knowing the Markov equation.  
For more information on the
classical Markov equation we refer to \cite{Cassels} and to the
original article of Markov \cite{Markov}. 
This equation also
appears in the work of Rudakov on vector bundles on $\Pp^2$
(\cite{Rudakov}) and has known applications  for mutations of
exceptional sequences (see e.~g.~\cite{GR}). The Markov
equation we study in this note differs slightly from the classical one: 
Markov constructed in \cite{Markov} all integral solutions of the equation
$$x^2 + y^2 + z^2 -3xyz = 0 $$
whereas we consider a rescaled version (dividing all variables by 3) and allow arbitrary right hand side $C$ (see section \ref{Sproofs}).

This paper is organized as follows. In Section \ref{Saction} we
define the group action we are interested in and prove some preliminary
lemmata. This closely follows \cite{Cassels}. In Section \ref{SMarkov},
we explain how to obtain the Markov equation and show its invariance
under cluster mutations. This is well known for mutations of exceptional
sequences (see \cite{Rudakov}) and the classical Markov equation.  We
consider the real points of the 
algebraic variety  defined by the Markov equation (it is a
differentiable manifold for $C\not=4$) and compute its connected
components in 
Section \ref{Sconnected}. It turns out that the group action respects
the components. 
Moreover, we need to compute the  
singularities. In Section \ref{Sfunddom} we construct a fundamental
domain for the action of the group $\oGamma$ on the set of all cluster-cyclic triples with real coefficients. Moreover, we show that the number of orbits for any
given Markov constant $C\not = 4$ is finite. Finally, in Section
\ref{Sproofs} we prove our main theorems and the corollary.

{\sc Aknowledgment. } This work was done during a visit of the third
author in Sherbrooke. He would like to take the opportunity to thank
the second author for his invitation and the representation theory
group for the nice working conditons.
 

\section{The group action and mutations}\label{Saction}

Let $\Gamma$ be the group which is freely generated by three generators
$\mu_1, \mu_2,$ and $\mu_3$ of order $2$

$$
\mu_1: (x,y,z) \mapsto (yz - x,y,z), \quad  \mu_2: (x,y,z) \mapsto (x,xz - y,z),\quad \mu_3: (x,y,z) \mapsto (x,y,xy - z)
$$ 
and let $\oGamma$ the semi-direct product of $\Gamma$ with the
symmetric group $S_3$ (see also below).
It follows from Lemma \ref{LMarkov1} that $\Gamma$ is the subgroup of $\oGamma$ which is
generated by all cluster mutations.
 Note that $\Gamma$ describes the cluster mutation as defined in \cite{FZ1} only when the quiver is cyclic.
In case of an acyclic quiver, the cluster mutation is given by  a different formula. 
Since we only want to decide
whether a cyclic quiver is cluster-cyclic or not, it is sufficient to study the
action of $\oGamma$ (or $\Gamma$) on $\Zz^3$ (or on $\Rr^3$). 
A triple $(x,y,z)$ is then cluster-acyclic precisely when there exists a
triple in the $\oGamma$-orbit $\oGamma (x,y,z)$ with one non-positive
entry.

 In case there are more than $3$ vertices, one could  also derive a
 version of the Markov equation in a manner similar to that described in
 Section \ref{SMarkov}. 
 However, one then has  to consider also
 acyclic quivers, and we see no 
way to make the group action and the Markov equation compatible with cluster mutations on acyclic quivers.


The group $\Gamma$ acts via cluster mutations on the three-dimensional affine
space $\Aa^3$ defined over any commutative ring $k$. Moreover,  the symmetric group $S_3$ acts on $\Gamma$ via
permutation of the generators: 
$$\sigma \cdot \mu_{a(l)}
\mu_{a(l-1)}\ldots \mu_{a(2)}\mu_{a(1)} :=
\mu_{\sigma(a(l))} \mu_{\sigma(a(l-1))}\dots
  \mu_{\sigma(a(2))}\mu_{\sigma(a(1))} \quad \mbox{ for } \sigma \in S_3 \;.
$$ 
Thus $\oGamma$
is the semi--direct product of $S_3$ with $\Gamma$. 
In Section \ref{Sfunddom}, we construct a fundamental
domain for the action of the group $\oGamma$ 
 on the set of all cluster-cyclic triples with Markov
constant less than $4$. 
For the cluster-cyclic triples  with $C=4$, however, we cannot construct a fundamental domain.
This is explained in detail in Section
\ref{Sconnected}, when we compute the connected
components and the singularities of the sets $V(C)$ over $\Rr$. 

In this section, we only prove some first elementary properties of the
group action. 
Note that the cluster mutation $\mu_i(x,y,z)$ changes only one of the three components $x,y$ and $z$.

So we can define $(x,y,z) \le
\mu_i(x,y,z)$ if all three components satisfy this inequality. It is a
partial order on the set of all real triples $(x,y,z)$. For us it is
important that we can compare $(x,y,z)$ with $\mu_i(x,y,z)$ for all
$i = 1,2,3$.

We distinguish the following
three cases: 

(M1) $(x,y,z) \leq \mu_i(x,y,z)$  for all $i = 1,2,3$, \\
(M2) $(x,y,z) \leq \mu_i(x,y,z)$  for precisely two indices $i$, \\
(M3) $(x,y,z) \leq \mu_i(x,y,z)$ for at most one index $i$.

Moreover, for any triple with $x \geq y \geq z$ the inequality $xy - z
< z$ already implies $xz - y < y$ and $xz -y < y$ implies $yz - x <
x$.

\begin{Lemma} \label{LMarkov1} 
Consider  three integers $x,y,z \geq 0$. \\
a) If $(x,y,z)  \not= (0,0,0)$ satisfies (M1), then $x,y,z \geq 2$ and
the triple is cluster-cyclic.
Moreover, in this case $(x,y,z)$ is the unique
triple in its $\oGamma$-orbit satisfying conditon (M1) whereas all other triples  in its $\oGamma-$orbit 
satisfy condition (M2). \\ 
b) In the $\oGamma-$orbit of $(x,y,z)$ there is either a unique triple satisfying (M1) or a triple
satisfying (M3). \\
c) If $(x,y,z)$ satisfies (M3), then $\min\{x,y,z\} < 2$, and the triple is
cluster-acyclic. \\
d) If $(x,y,z)$ is a triple where one entry
equals $2$, then the 
triple is either cluster-acyclic or it is of the form $(x,x,2)$ with
$x \geq 2$. \\
e) If $(x,y,z)$ is a fixed  point under the $\Gamma-$action,
then either \\
$C(x,y,z)=0$ and  $(x,y,z) = (0,0,0)$ or \\
$C(x,y,z)=4 $ and  $(x,y,z) = (2,2,2)$. \\
f) The only element $\gamma \in \oGamma$ acting trivially on
$\Aa^3(\Zz)$ is the unit element.
\end{Lemma}

{\sc Proof. }
Statement a) summarizes Markov's result from \cite{Markov}, see also
\cite{Cassels}. A proof adapted to the situation we consider here is
given in \cite{ABBS} from which also most of the other statements can
be derived. 
We give here the main arguments, starting with part c) of the lemma. Assume that  the triple
satisfies (M3) and suppose, without loss of generality, that $x \geq y
\geq z$.  
Then $yz - x < x$ and $xz - y < y$ (this holds up to the action of $S_3$, then we use the assumpion $x \geq y \geq z$ ) and thus $z^2 y < 2xz
< 4y$ and $z^2 < 4$. Consequently, $z = 0$ (the quiver is acyclic) or
$z = 1$. In the second case, we apply the mutation $(x,y,1)
\mapsto (y - x,y,1)$ to obtain an acyclic quiver. \\
To prove a), assume now we start with a triple $(x,y,z)$ satisfying
(M1). Thus $x,y,z \geq 2$, otherwise we can mutate the triple to a
smaller one. We start to mutate the triple as
long we get a triple satisfying (M2); the triples we obtain in each
step away from $(x,y,z)$ are then larger or equal to their
predecessors. 
Assume $(x',y',z')$ is the first
triple not satisfying (M2). Then it does not satisfy (M1), otherwise
the previous one does not satisfy (M2). It also can not satisfy (M3),
since otherwise (using part c)) we got a triple with one entry smaller
than $2$, but $(x',y',z') \geq (x,y,z)$. One can think about all triples as the vertices of a graph,
where each triples has at most three neighbours, the three mutated
triples. This graph is also ordered, that is, it defines a
poset. Using this poset it is not difficult to see the properties we
claimed.  \\
The uniqueness can also be shown as follows:
Assume there are two triples 
in the $\Gamma$--orbit
satisfying (M1). Then there exists a sequence of mutations from the
first to the second. Along  this sequence there must be at least one
triple which admits two smaller neighbours (use the poset as described
above), thus it satisfies (M3). But this a
contradiction to c). \\
To prove the first part of b) 
we note that we can always apply a mutation to obtain a smaller triple unless the triple 
 satisfies (M1). So, by applying a finite number of mutations, one arrives either at a triple satisfying (M1) or one obtains a triple
with some entry smaller than $2$, then some element in the orbit satisfies (M3). \\ 
To prove d) we consider the following sequence of mutations (where we
assume $x > y$): 
$$
(x,y,2) \mapsto (2y - x,y,2) \mapsto (2y - x,3y - 2x,2) \ldots 
$$
After a finite number of steps $ay - (a-1)x < 0$, so the triple is
cluster-acyclic. \\
If we have a fix point, then we are in case (M1), so
$z^2x = 4x$, $y^2 z = 4z$, and $x^2 y = 4y$. So all values must be $0$
or of absolute value $2$. Then one checks all possible values
explicitly. \\
Finally we prove f). There exist triples with $x > y > z$
satisfying 
(M1) (take for example $(5,4,3)$). Then for all elements $\gamma \in
\oGamma$, we find $\gamma (x,y,z) \not= (x,y,z)$; otherwise we find an 
element $\gamma'$ (that is a subsequence of the mutations defining
$\gamma$) with $\gamma' (x,y,z)$ satisfying (M3); we use the same argument
as in  part a). This yields a contradiction to the uniqueness in
part b).
\hfill $\Box$

\section{The Markov equation}\label{SMarkov}

In this section, we explain why the Markov constant is invariant under
mutations. First, we recall from \cite{Rudakov} and \cite{GR} the appearance of the Markov equation from the study of exceptional sequences of vector bundles or modules, where one obtains a
similar formula (for vector bundles on $\Pp^2$ the classical Markov
equation appears). 
Note that exceptional sequences admit a similar notion of mutation
together with an action of the braid group. It turns out that for
exceptional sequences of length 3 the semidirect product of those
mutations with the symmetric group is isomorphic to $\oGamma$. 

For exceptional mutations (of exceptional sequences of length $3$), one obtains
the Markov constant as 
follows. Consider the Cartan matrix $D$ of the endomorphism ring of an exceptional sequence $(E_1,E_2,E_3)$.
It has the form 
$$
D := 
\left(
\begin{array}{ccc}
1 & x & y \\
0 & 1 & z \\
0 & 0 & 1
\end{array}
\right) 
$$
with integral entries $x,y,z$.
Now compute its corresponding Coxeter matrix
$\Phi$, then the Markov constant is the
trace of $\Phi$ up to a constant:
$$
\Phi =  - D^t D^{-1} = 
- \left(
\begin{array}{ccc}
1 & 0 & 0 \\
x & 1 & 0 \\
y & z & 1
\end{array}
\right)
\left(
\begin{array}{ccc}
1 & -x & xz - y \\
0 & 1 & -z \\
0 & 0 & 1
\end{array}
\right) =
\left(
\begin{array}{ccc}
-1 & * & * \\
* & x^2 - 1 & * \\
* & * & y^2 + z^2 - xyz -1
\end{array}
\right)
$$
$$
\mTr(\Phi)  = x^2 + y^2 + z^2 - xyz - 3.
$$

It follows from tilting theory (see e.~g.~\cite{Ringel}) that the Coxeter matrices
of two
derived equivalent algebras are conjugate, so, in particular, their
traces are  equal.
Since the mutation of a full strongly exceptional sequence (it defines
a tilting module as the direct sum of the elements of the sequence)
yields an endomorphism ring which is derived equivalent, the Markov
constant must be invariant under such mutations. 
This argument gives an explanation for
why the Markov constant occurs; the invariance itself can be 
checked directly:

\begin{Lemma} \label{LMarkovinv}
The Markov constant $C(x,y,z) = x^2 + y^2 + z^2 - xyz$ is invariant
under all cluster
mutations $\mu_i$ for $i=1,2,$ and $3$ and invariant under the action
of $\oGamma$.
\end{Lemma}

{\sc Proof. }
Obviously, the Markov constant is invariant under any permutation of
the variables $x,y,z$. Then it is sufficient to check invariance under
the transformation $x \mapsto yz - x$:
$$
(yz - x)^2 + y^2 + z^2 - (yz - x)yz = x^2 + y^2 + z^2 - xyz. 
$$
\hfill $\Box$

From the Markov equation we get also the eigen values of the Coxter
transformation. We first note that there exists always a largest real
eigen value in the wild case (see \cite{delaPenaTakane}) or all eigen
values are roots of unity in the tame case. Moreover, the
characteristic polynomial 
is symmetric of degree three: $T^3 + (C(x,y,z) - 3) T^2 + (C(x,y,z) -
3) T + 1$ 
(see 
\cite{Lenzing} and note that the trace of $\Phi$ is $C(x,y,z) - 3
$). All together we obtain the following result.
 Note that the Markov constant already determines all eigen
values. Moreover, all the exceptions in the main theorems satisfy
$\lambda$ is not real ($\lambda$ is not real precisely when $0 <
C(x,y,z) < 4$).

\begin{Lemma}
The eigen values of the Coxeter transformation $\Phi$ are $-1,
\lambda, 1/\lambda$ for some element $\lambda$ satisfying either \\
1) $\lambda$ is real with $|\lambda| \geq 1$ and then $C(x,y,z) \geq
4$ or $C(x,y,z) \leq 0$ or \\
2) $\lambda$ is a complex number with $|\lambda | = 1$ and then $0
\leq C(x,y,z) \leq 4$. \\
Finally $-1 + \lambda + 1/\lambda = C(x,y,z) - 3$. 
\end{Lemma}

\begin{Lemma}\label{LfirstMarkovvalue}
a) Assume a triple $(x,y,z)$, with $x,y,z \geq 0$, is cluster-acyclic,
then $C(x,y,z) > 0$ or 
$C(x,y,z) = 0$ and $x,y,z = 0$. \\
b) Assume a triple $(x,y,z)$ is cluster-cyclic, then $C(x,y,z) < 4$ or
$C(x,y,z) = 4$ and the triple is in the $\oGamma$-orbit of $(u,u,2)$
for some $u \geq 2$.  
\end{Lemma}

{\sc Proof. }
Assume a triple is cluster-acyclic, then there exists an element
$(x,y,z)$ in its orbit with $x,y \geq 0$ and $z \leq 0$. Then,
$x^2\geq 0$, $y^2\geq 0$, $z^2\geq 0$, and $-xyz \geq 0$ and we get
equality only for 
$(0,0,0)$. Assume the triple is cluster-cyclic, then 
we can assume $x \geq y \geq z \geq 2$ and $yz -x \geq x$ (so that it is
the unique element satisfying (M1), it is also the unique minimal element in
the orbit, see Lemma 
\ref{LMarkov1}). We consider $z=2$. Then we obtain $2y \geq 2x$ 
so $x=y$ with $C(x,x,2) = 4$. If $z \geq 3$,  then $zy \geq 2x$. Now
we assume $x=y$,  then $C(x,y,z) \leq 2x^2 +9 - zx^2 \leq 0$ with
equality only for $x,y,z = 3$. If  $z \geq 3$ and $zy = 2x$, then
$C(x,y,z) = - ((\frac{z^2}{4} - 1)y^2 + z^2) < 0$. Finally note that
$V(C)_z$ (here we fix the value of $z$ and consider $V(C)_z$ as an
affine algebraic variety with coordinates $x,y$) is a hyperboloid,
thus convex and $C(x,y,z)$ is at most the 
value for $x=y$ or for $zy = 2x$ we computed above (all this is 
elementary computation, to see more details
compare with the results in Section \ref{Sconnected} and the
computation of the fundamental domain in Section \ref{Sfunddom}). 
\hfill $\Box$
\medskip

Using Lemma \ref{LfirstMarkovvalue} we have already proven Theorem
\ref{Tmain1} and Theorem \ref{Tmain2} for all values $C < 0$ and $C>
4$ (the details will be given later in Section \ref{Sproofs}). The
remaining cases need some further investigation, in 
particular, we need to determine the possible values of the Markov
constant between $0$ and $4$. Finally, note that $4$ is the maximal
possible value of the Markov constant for a cluster-cyclic triple. This
and Lemma \ref{LMarkov1} a) explains why we obtain the functions $m^+$
and $m^-$ as solutions of the Markov equation for $C=4$.

\section{Connected components}\label{Sconnected}

In this section we consider the Markov equation over the real numbers
and compute the slices 
$$
V(C)_z :=\{ (x,y) \in \Rr^2 \mid C(x,y,z) = C \}
$$
for a fixed value $z$. Then the Markov equation is quadratic in $x$
and $y$ and we consider the quadrics $V(C)_z$ for the various values
of $C$ and $z$. It turns out that $V(C)_z$ is an ellipsoid for $|z| <
2$,  a pair of lines for $|z| = 2$ and a hyperbola for $|z|
> 2$ (it might be empty for $|z| \leq 2$). Using the geometry of
$V(C)_z$ we can easily determine the 
connected components and the singularities of $V(C)$. It turns out
that 
$V(C)$ is smooth except for $C=4$ (over $\Cc$ it also has a
singularity for $C=0$ in $(0,0,0)$). If $C=4$ then it has the $4$ singularities
$(2,2,2), (2,-2,-2), (-2,2,-2),$ and $(-2,-2,2)$. Moreover, the 
$\Gamma$-action respects the connected components. Since we have for
$0 \leq C \leq 4$ five connected components in the smooth part $V(C)
\setminus V(C)_{\mathrm{sing}}$, one bounded and four unbounded, we
obtain both cluster-cyclic and cluster-acyclic orbits over
$\Rr$. We finally classify these orbits over the integers, and
 obtain all 
orbits for $0 \leq C \leq 4$. The remaining values of $C$ are easier
to handle (Lemma \ref{LfirstMarkovvalue}). We also note that it is
convenient to work 
with the $\Gamma$-action here since $S_3$ permutes some of the
components, however it preserves all components contained in
$\Rr^3_{\geq 0}$.

\begin{Lemma}\label{Lconncomp}
We fix $z$ and consider the restricted action induced by the
mutations $(x,y,z) \mapsto (yz-x,y,z)$ and $(x,y,z) \mapsto
(x,xz-y,z)$.\\
a) If $|z| < 2$ then $V(C)_z$ is an ellipsoid (hence
bounded) and for generic $x$, $y$ and $z$ the $\Gamma$-orbit is dense in
$V(C)_z$. In particular, there is no fundamental domain for the
restricted action. Moreover, the set $V(C)_z$ is empty precisely when
$C < z^2$. \\ 
b) If $|z| = 2$  then  
$V(C)_z$ consists of two lines for $C > 4$, of one line for $C=4$ and
is empty for $C<4$. If $C=4$ then the restricted action is trivial,
otherwise each generator interchanges the two components and the
composition of the two generators acts as a translation along the
lines. A fundamental domain for the restricted action on the set
$ \{ (x,y) \in \Rr^2 \mid x \not= y \}$ is 
$\{(x,y) \in \Rr^2 \setminus  (0,0) \mid y \leq 2x, y \leq 0 \}$. \\
c) If $|z| > 2$  then $V(C)_z$ is a hyperbola for $C \not= z^2$. It is
the union of two lines $x = \lambda y$ and $x = y/\lambda $ for $C =
z^2$, where $\lambda + 1/\lambda = z$. For $C < z^2$ one component
of $V(C)_z$ is contained in $\Rr^2_{> 0}$, the other one in $\Rr^2_{<
  0}$. It turns out that each orbit is unbounded in this case. For the
restricted action we have a fundamental domain.
\end{Lemma}

{\sc proof. }
Assume $|z| < 2$, thus $x^2 + y^2 -xyz \geq 0$ and  $C(x,y,z) = x^2 +
y^2 + z^2 -xyz \geq z^2$. To see that the $\Gamma$-orbit is dense
compute the eigenvalues of the linear transformation defined by the
generators. If they are not rational the orbit must be dense. 
 And since the coefficient of $xy$ is
smaller than $2$, the
equation defines an ellipsoid. This
shows a). To prove b) write the Markov equation as follows:
$C(x,y,z) = (x-y)^2 + 4 = C$. Thus $V(C)$ is empty for $C < 4$ and for
$C = 4$ we get 
one line, for $C>4$ we get two lines. The  action is
generated by $(x,y) \mapsto (2y - x,y)$ and $(x,y) \mapsto (x,2x - y)$, and their
composition is $(x,y) \mapsto (2y - x,3y - 2x)$. This is a linear
action preserving the orientation, so we only need to compute the
image $\Rr_{\geq 0} (-1,-2) $ of the line $\Rr_{\geq 0}(1,0) $ and a
fundamental domain is the cone generated by these two rays. To prove
c) we first compute the hyperboloids and the asymptotic
lines. Finally, we can also use the linearity of the action to get a
fundamental domain.
\hfill $\Box$
\medskip

Now we use the previous result to get information on the action of $\Gamma$ and
the components.

\begin{Thm}
Consider the action of $\Gamma$ on $V(C)$. \\
a) The set $V(C)$ is a smooth manifold for all $C \not= 4$ and for
$C=4$ it has four cone singularities. Denote by $V(C)^{{\mathrm
    smooth}}$ the open smooth part. \\
b) The group $\Gamma$ respects the connected components of $V(C)^{{\mathrm
    smooth}}$. \\
c) The number of connected components depends on the value of $C$ as follows:\\
\begin{center}
\begin{tabular}{c|c|c|c}
 & conn.~comp. & conn.~comp. & compact \\
 &  of $V(C)$ & ~of $V(C)^{{\mathrm smooth}}$ & conn.~comp. \\
\hline
$C  <0 $ & $4$ & $4$ & $0$\\
$C=0$  & $5$ & $5$ & $1$\\
$0 < C < 4$ & $5$ & $5$ & $1$\\
$C=4$ & $1$ & $5$ & $0$\\
$C> 4$ & $1$ & $1$ & $0$
\end{tabular}
\end{center}
d) For $C=0$ the compact connected component is the point
$(0,0,0)$. All compact connected components are contained in
$[-2,2]^3$. All orbits in a non-compact connected component have a
sequence of elements that converges to infinity. \\
e) All finite orbits are contained in the compact components or in the
singular part.
\end{Thm}

{\sc Proof. }
We compute the singularities over $\Cc$: consider the differential
$$
\partial C(x,y,z)/\partial x = 2x - yz.
$$
It must vanish for any choice of the coordinate, thus
$$
2x = yz, 2y = xz, \mbox{ and } 2z = xy \mbox{ yields } 4x =  y^2x, 4y
= z^2y \mbox{ and } 4z = x^2z.
$$
Then either $x=0$ and $z,y=0$ too, or $|y| = 2$ and $|x| = |z| =
2$, which yields the triples $(2,2,2), (2,-2,-2), (-2,2,-2),$ and $(-2,-2,2)$. 
Now we consider $V(C)$ over $\Rr$ again. Since
$(0,0,0)$ is an isolated point in 
$V(0)$ (consider the slices $V(C)_z$ we have already computed in the
previous lemma) it is not a singular point (over $\Rr$). This proves a) and
the first part of d). To
show b) we note that the group acts via differentiable (in fact even
algebraic) transformations, thus preserves the singularities. To
finish the proof of b) we need  to show c) first:  
Assume first $C \geq 4$. Then $V(C)_z$ is never empty and connected
for $|z| < 2$. Now vary $z$, then we obtain a path from a point in
$V(C)_z$ to each other slice $V(C)_{z'}$. Thus $V(C)$ is
connected. Assume now $0 < C < 4$. Consider $2 > |z| > 0$. For $z = 2 - \varepsilon$ and
$\varepsilon >0$ sufficiently small the
variety $V(C)_z$ is empty. On the other hand $V(C)_z$ is non-empty for
$z=0$ and $|z| > 2$. Playing this with all three coordinates we get
$5$ connected components, one is contained in $[-2,2]^3$ and is
compact, the four others are contained in $[2,\infty]^3$, \quad
$[-\infty,-2] \times [-\infty,-2] \times [2,\infty]$, $[-\infty,-2]
\times [2,\infty] \times [-\infty,-2]$, respectively $[2,\infty] \times
[-\infty,-2] \times [-\infty,-2]$.
 Note that no component can be
contained in $[-\infty,-2]\times [2,\infty] \times [2,\infty]$, since
the Markov constant is at least $4$. If $C < 0$ the compact
component vanishes and for $C=0$ it is just the origin. If $0 \leq C
\leq 4$ then we have both, a compact component and four unbounded
components. The arguments are similar to the ones above. One considers
the slices $V(C)_z$ for the various $z$ and applies the action of the
symmetric group and the possible sign changes. \\ 
Now we can
finish the proof of b): since the five components for $C< 4$
(respectively the components in the smooth part for $C=4$) respect the
inequality $|z|,|y|,|z| \leq 2$ in two variables, they must respect it
also in the third one by part c). Note that $\oGamma$ (or any
permutation) can change signs and permutes three of the non-compact
components, it only preserves the compact component and the one in
$[2, \infty]^3$.

\section{The fundamental domain} \label{Sfunddom}

The uniqueness in Lemma \ref{LMarkov1} a) suggests the existence of a
fundamental domain for the $\oGamma$-action on the set of
cluster-cyclic triples in $\Rr^3$. We note that the uniqueness was
proven only for integral triples. Indeed, there exist real triples
(they cannot be integral by Lemma \ref{LMarkov1})
that are cluster-cyclic and only have triples of type (M2) in their
orbit. We will show that those triples can only occur for the Markov
constant $4$. Thus for all real triples with fixed Markov constant $C$ there exists a unique triple  
satisfying (M1) only if the Markov constant is less than
4. Consequently,  there should exist a fundamental domain for
all cluster-acyclic triples with Markov constant strictly less than
$4$. We define
$$
F := \{ (x,y,z) \in \Rr^3 \mid x \geq y \geq z \geq 2, yz \geq 2x \}
\mbox{ and } 
$$
$$
F^{\circ}  := F \setminus \{(x,x,2) \mid x \geq 2\} = \{(x,y,z) \in  F \mid C(x,y,z) <4 \}. $$

Using Lemma \ref{LfirstMarkovvalue} we obtain $C(F) \subseteq
(-\infty,4]$ and $C(F^{\circ}) \subseteq (-\infty,4)$, since the
  elements in $F$ and $F^{\circ}$ are cluster-acyclic and satisfy (M1). 

\begin{Thm}\label{Tfunddom}
a) Let $(x,y,z)$ be a real cluster-cyclic triple. Then there exists a
unique triple $f(x,y,z)$ in $F$ that is in the closure (with respect to the
ordinary topology) of the orbit $\oGamma (x,y,z)$. If the triple $(x,y,z)$ is
integral or of Markov constant strictly less than $4$ then $f(x,y,z) \in
F \cap \oGamma (x,y,z).$\\
b) Consider the action of $\oGamma$ on the set $\{ (x,y,z) \mid x,y,z
> 2; \quad C(x,y,z) <4 \}$. Then $F^{\circ}$ is a fundamental domain for this
action and the set $S_3 \cdot F^{\circ}$ is a fundamental domain for the
$\Gamma$-action.
\end{Thm} 

{\sc Proof. }
We first note that Lemma \ref{LMarkov1} a), c), d), e), and f) are
also true for real triples. Only b) might be not. We obtain, by the
uniqueness, a 
fundamental domain for the action on the set of all triples which do
have a triple satisfying (M1) in their orbit. We show, it is already
the fundamental domain as claimed above. Obviously, each
element in $F$ satisfies (M1). We need to show that $\oGamma F^{\circ} = \{
(x,y,z) \mid x,y,z > 2; C(x,y,z) < 4 \}$. Take any triple in
$\{(x,y,z) \mid x,y,z > 2; C(x,y,z) \leq 4 \}$. Now we start to mutate
in the unique way, so that in 
each step we get a strictly smaller triple (using the partial order
defined in Section \ref{Saction}). Assume we got a smallest
triple $f(x,y,z)$, then it is an element in $F$ (since all elements
satisfying (M1) are in $F$ by definition). Assume there is not
a minimal triple with Markov constant $C < 4$. Then the sequence of
triples converges (it is decreasing and bounded  below) and the limit $f(x,y,z)$ must be in $F$. If the
triple $f(x,y,z)$ is in
the interior  of $F$ (this is $F$ without its boundary) then a small
ball around $f(x,y,z)$ is also in $F$ contradicting the fact that no
element of the convergent sequence is in $F$. A similar argument works
if $f(x,y,z)$ is in the boundary of $F$ for $C(x,y,z) = C(f(x,y,z)) <
4$, take a small  ball around the point $f(x,y,z)$ with all Markov
constants strictly less than $4$. Then we find that
the ball around $f(x,y,z)$ is contained in  $S_3 F \cup \mu_1 S_3F
\cup \mu_2 S_3 F \cup \mu_3 S_3 F$ (one can see, that each of the elements
$\mu_i$ fixes one of the boundaries. 
\hfill $\Box$

\begin{Lemma}
a) Assume $(x,y,z)$ is in $F_{\Zz}$. Then $C(x,y,z) = 4$, $C(x,y,z) =
0$ or $C(x,y,z) < 0$. \\
b) $F_{\Zz} \cap V(0) = \{(3,3,3) \} $\\
c) $F_{\Zz} \cap V(4) = \{(x,x,2) \mid x \in \Zz_{\geq 2}\} $.
\end{Lemma}

{\sc Proof. }
We need to compute the maximal value of the Markov constant on the
slices of the fundamental region $F_z := \{ (x,y,z) \in F \}$ where we
fix $z$ (we show that only $4$, $0$ and negative values may occure). Since $F_z$ is obviously bounded by lines, it is an
intersection of affine half spaces, in particular it is convex. It is
a little explicit computaton to see that the function $C(x,y,z)$ takes its minimal
value on the vertices of $F_z$, that is on $(z,z,z)$ or on
$(z^2/2,z,z)$.
The Markov constants are $3z^2 - z^3 = z^2(3 - z)$ and $- z^4/4 + 2z^2 = z^2(2 -
z^2/4)$.  Considering the values for $z = 2,3,4,\ldots$ we obtain that the
maximum is $z^2(3 - z)$ and it is in the critical part $[0,4]$ only for $z = 2$ or $3$. Finally, we need to classify the orbits
for $C=4$, that is already done using Lemma \ref{LfirstMarkovvalue}, and for
$C=0$. For $C=0$ it is the classical Markov equation and there are
precisely two orbits: $(0,0,0)$ (it is acyclic) and $\oGamma \cdot (3,3,3) =
\Gamma \cdot (3,3,3)$ (it is cyclic). \\
For convenience we repeat the arguments: assume $C(x,y,z) = 0$,
considering the equation modulo $3$ shows that $3$ divides each entry in
the triple. Now use the fundamental region. If $z \geq 4$  then there
is no integral point in $F$ with Markov constant $0$. We can assume
$z=3$, since $z = 2$ can not appear. Since $C = 0$ is the maximum on
$F_3$ all other elements in $F_{\geq 3}$ must have Markov constant
strictly smaller than $0$.
\hfill $\Box$

\begin{Cor}\label{Cfinnumborbits}
For each integer $C \not= 4$ the number of integral $\oGamma$-orbits
with Markov constant $C$ is finite. 
\end{Cor}

{\sc Proof.} The same arguments (convexity and shape of $F_z$) show
that the set $F_z \cap \{(x,y,z) \mid C(x,y,z) \geq C\}$ is bounded
for any $z$ and any $C$. Consequently, there are only finitely many
lattice points in this set and there are only finitely many orbits. If
the triple is cluster-acyclic we consider the set of all acyclic
triples with Markov constant C. This set is also finite, since the
equation
$$
x^2 + y^2 + z^2 - xyz = C
$$
has only finitely many integral solutions for $x,y \geq 0$ and $z \leq 0$.
\hfill $\Box$

\section{Proof of the Main Theorems}\label{Sproofs}

In this section we finally collect all the arguments for the proof of
the two main theorems. Since the Markov constant $C(x,y,z)$ is
invariant under cluster mutation and permutation we can consider the
various values of $C$ case by case. \\
Assume first $C>4$ or $C<0$. Then parts of the theorems follow from Lemma
\ref{LfirstMarkovvalue}: in Theorem \ref{Tmain1} (1), (2) and (3) are
equivalent. Since $m^+(x,y) = z$ or $m^-(x,y) =z$  precisely when
$C(x,y,z) = 4$ we obtain $m^+(x,y) < z$ or $m^-(x,y) > z$
precisely when $C(x,y,z) > 4$ for all $x,y \geq 2$. Thus, under the
assumption $C(x,y,z) > 4$ condition (1) is also equivalent to (4).
Similar arguments work for  Theorem \ref{Tmain2} under the assumption
$C(x,y,z) < 0$: we get (1) if and only if (2) and (3) by Lemma
\ref{LfirstMarkovvalue}. The same argument as above and $x,y,z \geq 2$
(otherwise $C(x,y,z) \geq 0$ by Lemma \ref{Lconncomp} a) and its
proof) shows the equivalence with (4). Finally (5) is equivalent by
Theorem \ref{Tfunddom} a). \\
Now we consider $C(x,y,z) = 0$. Then we are in the case of the
classical Markov equation
$$
X^2 + Y^2 + Z^2 -3XYZ = 0, \mbox{ where } 3X = x, 3Y = y, 3Z = z.
$$
We have already shown that each Markov triple is divisible by $3$, so
we can define new variables $X,Y,Z$ and get a new equation as
above. The solutions are well-known and the only orbits are $(0,0,0)$
(it is cluster-acyclic)
and $\oGamma (3,3,3) = \Gamma (3, 3,3)$ (it is
cluster-cyclic and corresponds to the orbit of (1,1,1) for the
classical Markov equation). Comparing with the statement in the
theorems finishes 
this case. \\
Now we consider $0 < C(x,y,z) < 4$. In this case there is no integral
triple in the fundamental domain $F$, thus all triples must be
cluster-acyclic. Using the connected components we can, case by case,
classify all triples with $x,y,z \leq 2$ explicitly (gives the list in
Theorem \ref{Tmain1} (3) and nothing in Theorem \ref{Tmain2}). \\
Finally, we need to consider $C(x,y,z) = 4$. Again we classify all
triples with $x,y,z \leq 2$ (we obtain the two orbits in Theorem \ref{Tmain1}
(3) and one orbit in Theorem \ref{Tmain2} (3).

\medskip
{{\small Andre Beineke \\
Fakult\"at f\"ur Mathematik \\
 Universit\"at Bielefeld \\
 POBox 100 131 \\
 D-33 501 Bielefeld \\
 Germany \\
 E-mail: abeineke@math.uni-bielefeld.de}}
\medskip

{{\small Thomas Br\"ustle\\
{ D\'epartement de
  Math\'ematiques, 
  Universit\'e de Sherbrooke\\
  Sherbrooke,  Qu\'ebec , J1K 2R1,  Canada  \quad and \\
   Department of Mathematics, Bishop's University\\
     2600 College St., Sherbrooke, Qu\'ebec, J1M 0C8, Canada} \\
 { E-mail:  thomas.brustle@usherbrooke.ca} \quad and \quad
 { tbruestl@ubishops.ca}}}
 \medskip

{{\small Lutz Hille}\\
{\small Freie Universit\"at Berlin}\\
{\small Fachbereich Mathematik und Informatik } \\
{\small II. Mathematisches Institut }\\
{\small Arnimallee 3}\\
{\small D-14195 Berlin }\\
{\small Germany}\\
{\small E-mail: hille@math.fu-berlin.de}}\\

%
\newpage

\centerline{\Large {Appendix}}

The aim of the appendix is to present an interpretation of the Markov constant
$C(x,y,z) = x^2+y^2+z^2-xyz$ in the mutation acyclic case
in terms of a first Hochschild cohomology group.

Let $H = K\mathcal{Q}$ be a basic hereditary algebra,
where $\mathcal Q$ is a finite  quiver without oriented
cycles, and $K$ is an algebraically closed field. The $n$
vertices of $\mathcal Q$ are $\{1,\ldots ,n\}$. The set of arrows
is denoted by $\mQ_1$. If $\alpha$ is an arrow, then $s(\alpha )$ denotes its
starting point and $t(\alpha )$ its terminating point. By $e_i$ we denote the 
primitive idempotent of $H$, corresponding to the vertex $i$.

Denote by $\Ho^1(H)\cong \Ext (H,H)$ the first Hochschild cohomology group of $H$
with coefficients in $H$, where $H^e$ is the envelopping algebra of $H$.
It is shown in \cite[1.6]{H} that $\dim \Ho^1(H)$ can be expressed as
$$\dim \Ho^1(H) = d-n +\sum_{\alpha\in \mQ_1} \nu (\alpha ),$$
where $d$ is the number of connected components of the quiver $\mQ$, and
$\nu (\alpha )$ is the number of paths from $s(\alpha )$ to $t(\alpha )$
which means $\nu (\alpha ) = \dim e_{t(\alpha )}He_{s(\alpha )}$.

It should be mentioned that, similar to the Markov constant, also $\dim \Ho^1(H)$
is related to the trace. of the Coxeter transformation of the hereditary algebra
$H$, see for example \cite[3.2.1]{L}.

Let $\mZ (x,y,z)$, with $x,y,z>0$ be a cyclic quiver with three vertices

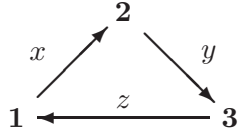
\begin{center}
\setlength{\unitlength}{2pt}
\begin{picture}(60,25)\thicklines
\put(0,0){\hin{$\bf{1}$}}
\put(20,20){\hin{$\bf{2}$}}
\put(40,0){\hin{$\bf{3}$}}
\put(20,3){\hin{$z$}}
\put(4,12){\hin{$x$}}
\put(36,12){\hin{$y$}}
\put(36,0){\vector(-1,0){32}}
\put(4,4){\vector(1,1){13}}
\put(24,16){\vector(1,-1){13}}
\end{picture}
\end{center}
\medskip

where $i \str{a}{\lra}j$ means that there are $a$ arrows from $i$ to $j$.
\medskip

If the  cyclic quiver $\mZ (x,y,z)$ is mutation acyclic, then, after a finite sequence of 
mutations one gets an acyclic quiver $\mQ (r,s,t)$ of the form
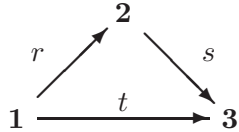
\begin{center}
\setlength{\unitlength}{2pt}
\begin{picture}(60,25)\thicklines
\put(0,0){\hin{$\bf{1}$}}
\put(20,20){\hin{$\bf{2}$}}
\put(40,0){\hin{$\bf{3}$}}
\put(20,3){\hin{$t$}}
\put(4,12){\hin{$r$}}
\put(36,12){\hin{$s$}}
\put(4,0){\vector(1,0){32}}
\put(4,4){\vector(1,1){13}}
\put(24,16){\vector(1,-1){13}}
\end{picture}
\end{center}
\medskip

with $r, s>0$ and $t\geq 0$. Hence $\mZ (x,y,z)$ is the quiver of a cluster tilted
algebra $\Gamma$  of type $H = K\mQ (r,s,t)$, see for example \cite{BMR}. In this case one has:
\medskip

\noindent{\bf Theorem } {\em If $\mZ (x,y,z)$ is the quiver of a cluster tilted algebra $\Gamma$
of type $H$,
where $H$ is connected hereditary of rank three, then
$$C(x,y,z) -2 = \dim \Ho^1(H).$$}
\medskip

\noindent{\em Proof. } Starting with the path algebra H of the quiver $\mQ (r,s,t)$ one gets from Happel's
result the formula 
$$ \dim \Ho^1(H) = r^2+s^2+t^2 +rst -2.$$
\medskip

Mutation at the vertex $2$ of the quiver $\mQ (r,s,t)$  yields a cyclic quiver $\mu_2 \mQ(r,s,t)$ of the
form

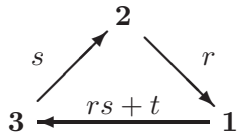
\begin{center}
\setlength{\unitlength}{2pt}
\begin{picture}(60,25)\thicklines
\put(0,0){\hin{$\bf{3}$}}
\put(20,20){\hin{$\bf{2}$}}
\put(40,0){\hin{$\bf{1}$}}
\put(20,3){\hin{$rs+t$}}
\put(4,12){\hin{$s$}}
\put(36,12){\hin{$r$}}
\put(36,0){\vector(-1,0){32}}
\put(4,4){\vector(1,1){13}}
\put(24,16){\vector(1,-1){13}}
\end{picture}
\end{center}
\medskip

But the Markov constant of this cyclic quiver is 
$$C(t+rs,s,r) = (t+rs)^2+s^2+r^2-(t+rs)rs =r^2+s^2+t^2+rst,$$
which proves the theorem, since the Markov constant is invariant under mutations, as long as the
quiver is cyclic.
\bigskip

\noindent{\bf Remark. } Let  $\mZ (x,y,z)$ be a fixed cyclic quiver, which is mutation acyclic, hence the quiver of a
cluster tilted algebra $\Gamma$. Using its Markov constant  $C(x,y,z)$ one can determine the finite list
of connected hereditary path-algebras $H$, with three vertices such that $\dim \Ho^1(H) +2 = C(x,y,z)$. $\Gamma$
then is cluster tilted of tpe $H$, where $H$ belongs to this list. The list can be quite big, as the
following example shows: Let $H$ be the path algebra of a quiver  of type $\mQ (2,s,t)$. Then
$\dim \Ho^1(H) = 2+(s+t)^2$.
%

{{\small Otto Kerner}\\
{\small Mathematisches Institut }\\
{\small Heinrich-Heine-Universit\"at D\"usseldorf } \\
{\small Universit\"atsstra\ss{}e 1   }\\
{\small D-40225 D\"usseldorf }\\
{\small Germany}}\\


\begin{thebibliography}{Dillo 50}

\bibitem{ABBS} I.~Assem, M.~Blais, T.~Br\"ustle, A.~Samson,
  \emph{Mutation classes of skew-symmetric $3 \times 3-$matrices},
  arXiv:math.RT/0610627, to appear in Comm.~Alg.
  
\bibitem{BMRRT} A. Buan, R. Marsh, M. Reineke, I. Reiten and
  G. Todorov, {\em Tilting theory and cluster combinatorics}, 
   Adv. Math.  {\bf 204} (2006), 572-612.

\bibitem{Cassels} J.~W.~S.~Cassels, {\emph An Introduction to
    Diophantine Approximation. } Facsimile reprint of the 1957
  edition. Cambridge Tracts in Mathematics and Mathematical Physics,
  No.~45.~Hafner Publishing Co., New York, 1972.

\bibitem{delaPenaTakane}J.~A.~de la Peña, M.~Takane, {\em Spectral
    properties of Coxeter transformations and applications. }  
Arch.~Math.~(Basel) 55 (1990), no.~2, 120--134.
 
\bibitem{FZ1} S. Fomin and A. Zelevinsky,
 {\em Cluster algebras I. Foundations}, J. Amer. Math. Soc. {\bf
 15}(2), (2002), 497-529 (electronic) 

\bibitem{FZ2}  S. Fomin and A. Zelevinsky, {\em Cluster algebras II. Finite type
  classification}, Inventiones Mathematicae {\bf 154}(1),
  (2003), 63-121. 


 
\bibitem{GR} A.~L.~Gorodentsev, A.~N.~Rudakov, \emph{ Exceptional
    vector bundles on projective spaces. }  Duke Math.~J.~54  (1987),
  no.~1, 115--130. 

\bibitem{K} O. Kerner, {\em Wild cluster tilted algebras of rank three}, preprint, July 2006.

\bibitem{Lenzing} H.~Lenzing, {\em
Coxeter transformations associated with finite-dimensional algebras. } Computational methods for representations of groups and algebras (Essen, 1997), 287--308,
Progr.~Math.~173, Birkhäuser, Basel, 1999. 

\bibitem{Markov} A. Markoff, \emph{Sur les formes quadratiques binaires ind\'efinies},  
Mathematische Annalen, Band XVII (1882), p. 379-399.

\bibitem{Ringel} C.~M.~Ringel, \emph{  Tame algebras and integral quadratic forms}. Lecture Notes in Mathematics, 1099. Springer-Verlag, Berlin, 1984. 

\bibitem{Rudakov} A.~N.~Rudakov
\emph{ Markov numbers and exceptional bundles on $P\sp 2$. } (Russian)
Izv.~Akad.~Nauk SSSR Ser.~Mat.~52 (1988), no.~1, 100--112, 240;
translation in Math.~USSR-Izv.~32 (1989), no.~1, 99--112. 


\end{thebibliography}

\begin{thebibliography}{ATTTS}
%

\bibitem{BMR}
{\sc A. Buan,R. Marsh and I. Reiten},
{\sl Cluster mutations via quiver representations},
arXiv:math.RT/0412077.
%
\bibitem{H}
{\sc D. Happel},
{\sl Hochschild cohomology of finite-dimensional algebras},
in Sem. M.P. Malliavin, Springer Lect. Notes Math. {\bf 1404} (1989), 108--126.
%
\bibitem{L}
{\sc F. Lukas},
{\sl Elementare Moduln \"uber wilden erblichen Algebren},
Thesis, D\"usseldorf 1992.
%
\end{thebibliography}
\end{document}